\documentclass[final,5p,twocolumn]{elsarticle}

\usepackage{graphicx}          
\usepackage[dvips]{epsfig}    

\usepackage{amssymb}
\usepackage{amsthm}

\usepackage[table]{xcolor}

\usepackage{amsmath}          

\usepackage{algorithmic}

\journal{ArXiv}

\begin{document}

\begin{frontmatter}

\title{Energy-saving sub-optimal sliding mode control with bounded
actuation}

\author[label1]{\scriptsize Michael Ruderman}
\ead{michael.ruderman@uia.no}
\author[label2]{\scriptsize Alessandro Pisano}
\author[label2]{\scriptsize Elio Usai}

\address[label1]{\scriptsize Department of Engineering Sciences, University of Agder, Norway}
\address[label2]{\scriptsize Department of Electrical and Electronic Engineering, University of Cagliari, Italy}

\normalsize

\begin{abstract}                          
The second-order \emph{sub-optimal} sliding mode control (SMC),
known in the literature for the last two decades, is extended by a
control-off mode which allows for saving energy during the finite
time convergence. The systems with relative degree two between the
sliding variable and switching control with bounded actuation are
considered, while the matched upper-bounded perturbations are not
necessarily continuous. Detailed analysis of the proposed
energy-saving sub-optimal SMC is performed with regard to the
parametric conditions, reaching and convergence time, and residual
steady oscillations if the parasitic actuator dynamics is added.
Constraints for both switching threshold parameters are formulated
with respect to the control authority and perturbations upper
bound. Based on the estimated finite convergence time, the
parameterization of the switching thresholds is solved as
constrained minimization of the derived energy cost function. The
total energy consuming control-on time is guaranteed to be lower
than the upper-bounded convergence time of the conventional
sub-optimal SMC. Numerical evaluations expose the properties of
the proposed energy-saving sub-optimal SMC and compare it with
conventional sub-optimal SMC in terms of the fuel consumption
during the convergence.
\end{abstract}


\end{frontmatter}

\newtheorem{thm}{Theorem}
\newtheorem{lem}[thm]{Lemma}
\newtheorem{clr}{Corollary}
\newdefinition{rmk}{Remark}
\newproof{pf}{Proof}

\section{Introduction}
\label{sec:1}

Sliding mode control (SMC), see \cite{utkin1992},
\cite{edwards1998}, \cite{shtessel2014} for fundamentals, is one
of the most promising robust control techniques. In SMC, the
control authority aims to force an uncertain system dynamics onto
a specified manifold $\sigma$ and then maintain the system
behavior on it, this way ensuring a reduced-order dynamics and
insensitivity to a certain class of perturbations. Being in
sliding mode implies $\sigma=0$ for all times $t > T_c$ after the
convergence at $0 < T_c < \infty$. In second-order sliding modes
\cite{levant1993}, \cite{fridman2002}, also the time derivative of
the sliding variable is forced to zero, i.e.
$\sigma=\dot{\sigma}=0$, thus allowing for robust control of
uncertain systems with relative degree two between the sliding
variable and control signal. Worth recalling is that a
second-order SMC belongs to a more generic class of higher-order
sliding mode approaches, which have been actively developed over
the last two decades for uncertain dynamic systems with relative
degree two and higher, see e.g. \cite{shtessel2014} and references
therein.

Among the second-order SMC issues, two of them can be highlighted
as particularly relevant for several applications. First, the time
derivative of the sliding variable may be hardy available in the
systems under control. Second, a discontinuous SMC implies a
high-frequency (theoretically infinite frequency) commutation of
the control variable, that can be both hardware-fatiguing and
energy-inefficient for different system plants. Addressing the
first above mentioned issue, the so called sub-optimal SMC
\cite{bartolini1997}, \cite{bartolini2003}, which is in focus of
our present work, was developed based on the bang-bang control
principles. Recall that the single commutation in an unperturbed
system characterizes the optimal bang-bang control. In turn, the
sub-optimal SMC converges robustly and in finite time to origin of
the $(\sigma,\dot{\sigma})$ plane by executing a commutating
control sequence $u(t)$ with increasing frequency. The
corresponding cumulative energy of a sub-optimal SMC is
proportional to the convergence time, since a discontinuous
commutating sequence implies the control signal is always on, i.e.
$|u(t)| \neq 0$ for $0 < t \leq T_c$. Here we explicitly point out
that the control behavior after the convergence, i.e. for $t
> T_c$, is not the subject of our recent discussion, while a
related remark is given at the end of the article. To the best of
our knowledge, an energy saving approach which allows also for
$u=0$ mode during convergence was not elaborated for second-order
SMC. At the same time, one should recall that an equivalent
fuel-optimal problem is well known for unperturbed second-order
systems \cite{athans1964}.

Against the above background, and being motivated by the previous
works \cite{athans1964}, \cite{bartolini1999}, we propose an
extension of the sub-optimal SMC which allows for energy saving
during the entire convergence phase. Here the systems with
relative degree two and discrete bounded actuation are targeted.
The latter implies the control value to be in the finite set $u
\in \{-U,\, 0,\, U \}$. The original sub-optimal SMC is modified
by introducing an additional switching threshold which enables for
a mandatory $u=0$ phase between two consecutive extreme values of
the sliding variable, cf. with \cite{bartolini2003}. The proposed
extension in the control law is rather simple, while the
associated analysis developed in detail is essential for
parameterization and guaranteed energy-saving operation of SMC.

The paper is organized as follows. The preliminaries of the
sub-optimal SMC are given in Section \ref{sec:2}, in close
accordance with its basic notations from \cite{bartolini2003}.
This section provides also the problem under consideration and the
formulated requirements for an energy-saving modification of the
sub-optimal SMC to be introduced. Section \ref{sec:3} contains the
main results. The proposed energy-saving sub-optimal SMC is
formulated and analyzed in terms of the convergence conditions and
a rigorous estimation of the worst-case reaching- and
convergence-time. Based on that, the constrained optimization
problem of energy-saving parametrization is formulated, and the
switching threshold parameters are determined so as to minimize
the fuel consumption. Also the well-known problem of residual
steady-state oscillations (i.e. chattering) owing to the
additional parasitic actuator dynamics is addressed in Section
\ref{sec:4}. This is done by applying the describing function
analysis of harmonic balance. Detailed numerical examples,
revealing and visualizing the main properties of the energy saving
sub-optimal SMC and also comparing it to the original sub-optimal
SMC, can be found in section \ref{sec:5}. Finally, section
\ref{sec:6} provides main concluding remarks.

\section{Preliminaries}
\label{sec:2}

Considered is a class of uncertain dynamic systems with the
matched bounded perturbations, for which the well-defined sliding
variable
\begin{equation}\label{eq:1}
\sigma(t) = \sigma \bigl( \mathbf{x}(t), t \bigr)
\end{equation}
has the relative degree $r=2$ with respect to the control variable
$u \in \mathbb{R}$. The measurable vector of the system states is
$\mathbf{x} \in \mathbb{R}^n$, with $n \in \mathbb{N}$ and $2 \leq
n < \infty$. The dynamics of the sliding variable can be written
in the normal form as
\begin{equation}\label{eq:2}
\ddot{\sigma}(t) = f(\cdot,t) + g(\cdot,t) u(t),
\end{equation}
where $f$ is a perturbation function and $g$ is the uncertain
input coupling function. Both are satisfying the global
boundedness conditions
\begin{itemize}
    \item[(i)] $\bigl | f(\cdot,t)  \bigr | \leq \Phi$, with $\Phi > 0$
    is a known real constant;
    \item[(ii)] $0 < \Gamma_m \leq  g(\cdot,t) \leq  \Gamma_M$, with $\Gamma_m, \Gamma_M$ to be known.
\end{itemize}
The sub-optimal second-order sliding mode control (SMC), proposed
initially in \cite{bartolini1997}, can be written in the following
form, cf. \cite{bartolini2003}
\begin{eqnarray}
\label{eq:3}
  u(t) &=& -\alpha(t) U \mathrm{sign}(\sigma - \beta \sigma_M),
  \\[1mm]
\label{eq:3b}
  \alpha(t) &=& \left\{%
\begin{array}{ll}
    1,        & \hbox{ if } (\sigma - \beta \sigma_M) \sigma_M \geq 0
    \\[1mm]
    \alpha^*, & \hbox{ if } (\sigma - \beta \sigma_M) \sigma_M < 0 , \\
\end{array}%
\right.    \\
\label{eq:3c}
  \beta & \in & [0;\, 1).
\end{eqnarray}
Here $U > 0$ is the minimal control magnitude and $\alpha^* \geq
1$ and $\beta$ are the modulation and anticipation factors,
correspondingly. The dynamic state $\sigma_M$ constitutes the last
\emph{extreme value} of the sliding variable during the control
\eqref{eq:3} operates on \eqref{eq:2}. The extreme value refers to
the value of $\sigma$ at the last time instant at which a local
maximum or horizontal flex point of $\sigma(t)$ has occurred, cf.
\cite{bartolini2003}. Following the previous works on the
sub-optimal SMC, see in survey \cite{bartolini2003}, the control
parameters must be set so as to satisfy
\begin{eqnarray}
\label{eq:4}
  U & > & \frac{\Phi}{\Gamma_m},\\[1mm]
  \alpha^* & \in & [1; \infty) \cap \Biggl ( \frac{2\Phi + (1-\beta)\Gamma_M U}{(1+\beta)\Gamma_m U}; \infty
  \Biggr).
\label{eq:5}
\end{eqnarray}
The condition \eqref{eq:4} represents the \emph{control
authority}, which is required so as to overcome the unknown
perturbations bounded by $\Phi$ and, therefore, to decide the sign
of $\ddot{\sigma}$. The condition \eqref{eq:5} determines and
ensures the convergence of $\sigma(t)$ to the origin, where the
second-order sliding mode appears in its proper sense. In
addition, a stronger inequality than \eqref{eq:5} can be imposed
by
\begin{equation}\label{eq:6}
\alpha^* \in  [1; \infty) \cap \Biggl ( \frac{\Phi +
(1-\beta)\Gamma_M U}{\beta\Gamma_m U}; \infty \Biggr),
\end{equation}
which is required for a monotonic convergence to zero, cf.
\cite{bartolini2003}. The latter means $\sigma(t)$ will experience
at most one zero-crossing, depending on the initial conditions
$(\sigma, \dot{\sigma})(0)$ before the first $\sigma_M$-state
occurs. For the sake of better distinguishing we will denote the
parametric conditions \eqref{eq:4}, \eqref{eq:5} by the
\emph{twisting convergence}, and \eqref{eq:4}, \eqref{eq:6} by the
\emph{monotonic convergence} conditions. Either of both conditions
ensure the establishment of the second-order sliding mode of
\eqref{eq:2} with \eqref{eq:3} in a finite time.

Worth recalling is also that $\dot{\sigma}$ is not available, so
that it is impossible to determine $\sigma_M$ by observing
zero-crossing of the $\dot{\sigma}$ state variable.
Notwithstanding, according to the 'Algorithm 2' provided in
\cite{bartolini1997}, a recursive detection of the extreme values
of $\sigma(t)$ takes place simultaneously with execution of the
control algorithm \eqref{eq:3}.

\subsection*{Problem under consideration}
\label{sec:2:sub:1}

The problem under consideration requires a bounded actuation $|u|
\leq U$ with only three discrete values, i.e. $u \in \{-U,\, 0,\,
U\}$, where $U > 0$ is a known parameter of the control system.
For such systems, we relax the above uncertainty condition (ii)
and assume a constant input coupling, meaning $\Gamma_m = \Gamma_M
= g = 1$. This implies the maximal control magnitude $U$ and leads
to $\alpha = 1$ to be substituted instead of \eqref{eq:3b}. Under
the above assumptions, the inequalities \eqref{eq:4}, \eqref{eq:5}
and \eqref{eq:6} result in a set of the parametric conditions
\begin{eqnarray}
\label{eq:7}
U & > & \Phi, \\
\label{eq:8}
\beta & > &  \frac{\Phi}{U} \qquad \qquad \hbox{: twisting convergence}, \\
\beta & > &  \frac{\Phi + U}{2U} \qquad \hbox{: monotonic
convergence}. \label{eq:9}
\end{eqnarray}
Now, we are in the position to specify an \emph{energy-saving
sub-optimal} SMC with the bounded actuation.

Given the fact of a bounded control authority $U$, we are
interested to modify the control \eqref{eq:3} with $\alpha=1$ so
as to allow for energy-saving control phases with $u=0$. Within
each reaching cycle, i.e. between two successive extreme values
$\sigma_{M_i}$ and $\sigma_{M_{i+1}}$ with $i \in \mathbb{N} <
\infty$, the new control policy must admit the following. The
control is on, i.e. $|u| = U$, on the intervals $\sigma \in
[\beta_1 \sigma_{M_i}; \sigma_{M_i}]$ and $\sigma \in
[\sigma_{M_{i+1}}; \beta_2 \sigma_{M_i}]$. And the control is off,
i.e. $u=0$, on the interval $\sigma \in [\beta_2 \sigma_{M_i};
\beta_1 \sigma_{M_i}]$, while $\beta_2 \neq \beta_1$ is required
for an energy-saving operation. The case where $\beta_2 = \beta_1
= \beta$ should correspond to the sub-optimal SMC \eqref{eq:3}
with $\alpha=1$, thus satisfying \eqref{eq:7} and either
\eqref{eq:8} or \eqref{eq:9}. This benchmarking control will be
denoted as \emph{conventional sub-optimal} SMC. For the
energy-saving sub-optimal SMC fulfills its task, following must be
guaranteed:
\begin{itemize}
    \item[(a)] convergence in finite time $T_c < \infty$, which implies $|\sigma_{M_{i+1}}| < |\sigma_{M_{i}}|$ for all $1 \leq i <
    \infty$;
    \item[(b)] fuel- and thus energy-saving, which implies
    $$
    \int \limits_0^{T_c} \bigl|u(\beta_1,\beta_2)\bigr| dt < \int \limits_0^{\hat{T}_c} \bigl|\hat{u}(\beta_1)\bigr|
    dt,
    $$
    where the variables with hat (like $\hat{\circ}$) are used for the
    conventional sub-optimal SMC.
\end{itemize}

\section{Energy-saving sub-optimal SMC}
\label{sec:3}

\subsection{Proposed control algorithm}
\label{sec:3:sub:1}

The proposed modification of the second-order sub-optimal SMC
algorithm is
\begin{equation}\label{eq:9a}
u(t) = -0.5 U \mathrm{sign}(\sigma - \beta_1 \sigma_M) - 0.5 U
\mathrm{sign}(\sigma - \beta_2 \sigma_M).
\end{equation}
In addition to \eqref{eq:7}, the thresholds relationship $\beta_1
> \beta_2$ is required, while the convergence and energy-saving
conditions on $(\beta_1,\beta_2)$ parameters are derived and
analyzed in the following. An initializing control action, cf.
\cite{bartolini2003},
\begin{equation}\label{eq:9b}
\bar{u}(t) = - U \mathrm{sign} \bigl(\sigma(t) - \sigma(0) \bigr),
\quad \forall \, t \in [0;t_{M_1}]
\end{equation}
is also required for speeding up the reaching of the first extreme
value at $t_{M_1}$.

\subsection{Proof of convergence and parametric conditions}
\label{sec:3:sub:2}

\begin{figure}[!h]
\centering
\includegraphics[width=0.99\columnwidth]{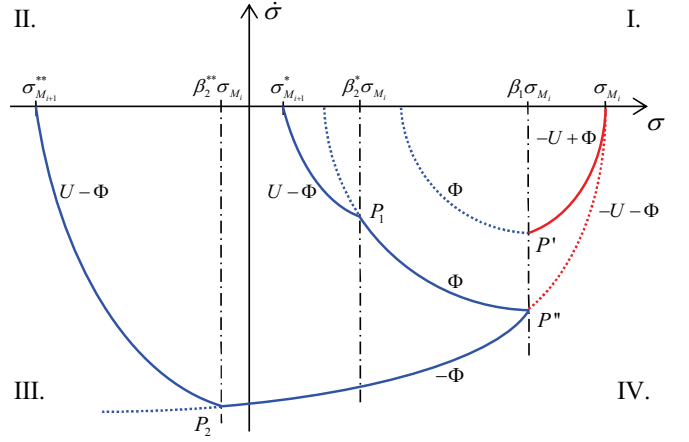}
\caption{Phase-plane of the energy-saving sub-optimal SMC for
worst-case scenario of the reaching time.} \label{fig:22}
\end{figure}

Using the general formula
\begin{equation}\label{eq:10}
\sigma = \sigma_0 + 0.5\, \frac{ \dot{\sigma}^2 -
\dot{\sigma}^2_0}{ v},
\end{equation}
for a double-integrator with $v = \pm U \pm \Phi$, the $(\sigma,
\dot{\sigma})$-trajectories can be calculated for an arbitrary
initial state $(\sigma_0, \dot{\sigma}_0)$. Then, without loss of
generality, consider the first extreme value $\sigma_{M_i} > 0$,
while the preceding reaching phase is always ensured by the
initialization control \eqref{eq:9b}. After applying first the
control value $u=-U$, according to \eqref{eq:9a}, the state
trajectory proceeds in the IV-th quadrant, cf. Figure
\ref{fig:22}. While both limiting cases of the perturbation $f =
\pm \Phi$ are effective, we focus here on the worst-case scenario,
that is essential for convergence, in which the sign of $f$ is
negative. Note that this renders the system to be underdamped.
Using $-U - \Phi$ in the denominator of \eqref{eq:10} and
substituting the threshold value $\beta_1 \sigma_{M_i}$ on the
left-hand side of \eqref{eq:10} one obtains the first peaking
point (labeled by $P''$) of the $\dot{\sigma}$-state as
\begin{equation}\label{eq:11}
\dot{\sigma}_{P''} = - \sqrt{2 \sigma_{M_i} (U+\Phi)(1-\beta_1)},
\end{equation}
At this point, the control is switching from $-U$ to $0$. Then,
evaluating the parabolic trajectory \eqref{eq:10} between the
$\beta_1 \sigma_{M_i}$ and $\beta_2 \sigma_{M_i}$ points, with
$u=0$ and worst-case $f=-\Phi$, one obtains
\begin{equation}\label{eq:12}
\beta_2 \sigma_{M_i} = \beta_1 \sigma_{M_i} + 0.5\, \frac{
\dot{\sigma}^2 - \dot{\sigma}_{P''}^2}{ -\Phi}.
\end{equation}
Substituting \eqref{eq:11} into \eqref{eq:12} and solving
\eqref{eq:12} with respect to $\dot{\sigma}$ results in
\begin{equation}\label{eq:13}
\dot{\sigma}_{P_2} = - \sqrt{ 2 \sigma_{M_i} \bigl( U(1-\beta_1) +
\Phi(1-\beta_2) \bigr) },
\end{equation}
where the second peaking point $P_2$ is the one where the control
is switching back from $0$ to $U$. Continuing this line of
calculations, one considers also the parabolic trajectory
\eqref{eq:10} after passing the $\beta_2 \sigma_{M_i}$-threshold
as
\begin{equation}\label{eq:14}
\sigma = \beta_2 \sigma_{M_i} + 0.5\, \frac{ \dot{\sigma}^2 -
\dot{\sigma}_{P_2}^2}{U -\Phi}.
\end{equation}
Here again, $f=-\Phi$ is assumed as worst-case from a convergence
point of view. Substituting the left-hand-side of \eqref{eq:14} by
the next extreme value $\sigma_{M_{i+1}}$, at which
$\dot{\sigma}=0$, results in
\begin{equation}\label{eq:15}
\sigma_{M_{i+1}} = \beta_2 \sigma_{M_i} + 0.5\,
\frac{-\dot{\sigma}_{P_2}^2}{U - \Phi}.
\end{equation}
Note that for convergence, it is strictly required that
$$
|\sigma_{M_{i+1}}| < |\sigma_{M_{i}}|
$$
for all pairs of the consecutive extreme values, indexed by $i$
and $i+1$. Since $\sigma_{M_{i}} > \sigma_{M_{i+1}}$ is guaranteed
for all trajectories in the IV-th and III-rd quadrants, cf. Figure
\ref{fig:22}, the above convergence inequality can be written,
with respect to \eqref{eq:15}, as
\begin{equation}\label{eq:16}
\beta_2 \sigma_{M_i} - 0.5\, \frac{\dot{\sigma}_{P_2}^2}{U - \Phi}
> - \sigma_{M_i}.
\end{equation}
Substituting \eqref{eq:13} into \eqref{eq:16} yields
\begin{equation}\label{eq:17}
1+ \beta_2 > \frac{\Phi + U - \Phi \beta_2 - U \beta_1}{U - \Phi},
\end{equation}
which results in the strict parameters condition
\begin{equation}\label{eq:18}
\beta_1 + \beta_2 > \frac{2 \Phi}{U}.
\end{equation}
Furthermore, one needs to ensure
\begin{eqnarray}
\label{eq:19}
  0 & \leq \beta_1 < & 1, \\
  -1 & < \beta_2 < & \beta_1,
\label{eq:20}
\end{eqnarray}
in order to realize the both control switchings between two
successive extreme values. Note that \eqref{eq:19} is inherited
from the conventional sub-optimal SMC, cf. \eqref{eq:3c}. Here we
emphasize that the parametric conditions \eqref{eq:18},
\eqref{eq:19}, \eqref{eq:20} must be satisfied for the
energy-saving sub-optimal SMC to yield the global convergence, for
$U > \Phi$ and irrespective the perturbations $|f(t)| < \Phi$.

\begin{figure}[!h]
\centering
\includegraphics[width=0.6\columnwidth]{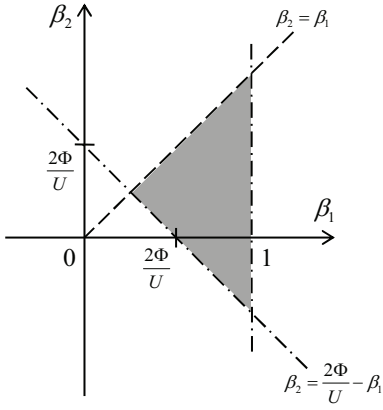}
\caption{Parametric constraints \eqref{eq:18}, \eqref{eq:19},
\eqref{eq:20} in $(\beta_1,\beta_2)$ plane.} \label{fig:21}
\end{figure}

The imposed $\beta_1, \beta_2$ constraints \eqref{eq:18},
\eqref{eq:19}, \eqref{eq:20} have an illustrative graphical
interpretation as shown in Figure \ref{fig:21}.
The $\beta_1, \beta_2$ values admissible for convergence are
inside of a grey-shadowed triangle, while on the $\beta_1 =
\beta_2$ edge the energy-saving control \eqref{eq:9a},
\eqref{eq:9b} reduces to the conventional sub-optimal SMC.
Obviously, the $\Phi/U$ ratio determines the size of an admissible
$\{\beta_1, \beta_2\}$ set, including whether the negative
$\beta_2$-values are allowed. It is also worth noting that the
opposite $\beta_2=2\Phi/U - \beta_1$ edge maximizes the
$\beta_1-\beta_2$ distance. However, this can largely increase the
convergence time (including $T_c \rightarrow \infty$), similar as
in case of a minimal-fuel control with free response time of an
unperturbed double-integrator analyzed in detail in
\cite{athans1964}.

\subsection{Reaching and convergence time}
\label{sec:3:sub:3}

Consider first the reaching time between two consecutive extreme
values $\sigma_{M_i}$ and $\sigma_{M_{i+1}}$. Recall that a
conservative approach is taken so that a worst-case scenario is
calculated for all control phases, in terms of a slowest time. The
associated trajectories of one reaching phase are schematically
drawn in Figure \ref{fig:22}. Similar as before, and without loss
of generality, consider the first extreme value $\sigma_{M_i} > 0$
so that the trajectories proceed in the III-rd and IV-th
quadrants. We recall that the threshold value $\beta_2$ can be
either positive or negative, cf. Figure \ref{fig:21}. Evaluating
first the peaking point $P'$ one obtains the corresponding
$\dot{\sigma}$-value as
\begin{equation}\label{eq:51}
\dot{\sigma}_{P'} = - \sqrt{2 \sigma_{M_i} (U-\Phi)(1-\beta_1)},
\end{equation}
which is reached by the slowest trajectory driven by $-U+\Phi$.
For the subsequent phase with control-off, i.e. $u=0$, one needs
distinguishing between two boundary trajectories, one driven by
$\Phi$ and another driven by $-\Phi$. Also we notice that both
control-off phases are considered starting from the same $P''$
peaking point. Indeed, one can recognize that starting from
$\beta_1 \sigma_{M_i}$, the $\Phi$-driven trajectory reaches the
$\sigma$-abscissa later from $P''$ than from $P'$ peaking point,
see Figure \ref{fig:22}. For the control-off phase between the
$\beta_1 \sigma_{M_i}$ and $\beta_2 \sigma_{M_i}$ thresholds,
starting from $P''$ and ending in $P_1$ peaking point, one obtains
\begin{equation}\label{eq:521}
\dot{\sigma}_{P_1} = -\sqrt{2 \sigma_{M_i} \bigl( U(1-\beta_1) +
\Phi(1-2\beta_1 + \beta_2) \bigr)}.
\end{equation}
Note that the second term under the square root in \eqref{eq:521}
must be positive, so here $\beta_2$ has to be larger than the
point where the $\Phi$-driven trajectory reaches the
$\sigma$-abscissa (see the dotted arc from $P_1$ in Figure
\ref{fig:22}).

With the determined $\dot{\sigma}$ values at all relevant peaking
points, i.e. $P'$, $P''$, $P_1$ and $P_2$, one can analyze the
slowest time for each phase of a reaching cycle individually.

The slowest time of the first control-on phase is
\begin{equation}\label{eq:522}
T' \equiv t_{P'} - t_{\sigma_{M_i}} =
\frac{|\dot{\sigma}_{P'}|}{U-\Phi} = \frac{\sqrt{2 \sigma_{M_i}} }
{U-\Phi} \sqrt{(U-\Phi)(1-\beta_1)}.
\end{equation}
In a similar way, one obtains the slowest time of the subsequent
control-off phase driven by $\Phi$ as
\begin{equation}\label{eq:523}
\begin{split}
T_1^* & \equiv t_{P_1} - t_{P''} =
\frac{|\dot{\sigma}_{P''}-\dot{\sigma}_{P_1}|}{\Phi} =
\frac{\sqrt{2 \sigma_{M_i}} } {\Phi} \, \times \\
& \Bigl( \sqrt{(U+\Phi)(1-\beta_1)}  - \sqrt{ U(1-\beta_1) +
\Phi(1-2\beta_1 + \beta_2)}  \Bigr).
\end{split}
\end{equation}
One can recognize that if $\beta_2=\beta_1$ (yielding the
conventional sub-optimal SMC), then both square root terms in
\eqref{eq:523} become equal, and $T_1^*$ becomes zero. One can
also recognize that the upper bound of $T_1^*$ is when the second
square root term in \eqref{eq:523} becomes zero. This results in
the maximal control-off time of this phase equal to
\begin{equation}\label{eq:524}
T^*_{1,\max} = \frac{\sqrt{2 \sigma_{M_i}} } {\Phi}
\sqrt{(U+\Phi)(1-\beta_1)}.
\end{equation}
When evaluating the control-off phase driven by $-\Phi$, cf.
Figure \ref{fig:22}, we write
\begin{equation}\label{eq:525}
\dot{\sigma}_{P_2} = \dot{\sigma}_{P''} - \Phi T^{*}_2.
\end{equation}
Substituting \eqref{eq:11} and \eqref{eq:13} into \eqref{eq:525}
and solving it with respect to $T^{*}_2$ results in
\begin{equation}\label{eq:526}
\begin{split}
T^{*}_2 = & \frac{\sqrt{2 \sigma_{M_i}} } {\Phi}  \, \times \\
& \Bigl ( \sqrt{U(1-\beta_1) + \Phi(1-\beta_2)} -
\sqrt{(U+\Phi)(1-\beta_1)} \Bigr). \end{split}
\end{equation}
Here, it can also be recognized that $T^{*}_2$ has its maximal
value when the second square root term in \eqref{eq:526} becomes
zero. Thus, the maximal control-off time for this phase is
\begin{equation}\label{eq:5266}
T^{*}_{2,\max} = \frac{\sqrt{2 \sigma_{M_i}} } {\Phi}
\sqrt{U(1-\beta_1) + \Phi(1-\beta_2)}.
\end{equation}
Comparing both maximal time values, one can show that
$$
T^*_{1,\max} < T^*_{2,\max}
$$
holds always due to $\beta_2 < \beta_1$. Despite this fact, which
determines the slowest control-off phase during one reaching
cycle, we recall that the overall convergence time depends also on
the number of cycles. Therefore, both \eqref{eq:524} and
\eqref{eq:5266} are taken into account in further calculations.

Evaluating both successive control-on phases, starting once from
$P_1$ and once from $P_2$ peaking point, we obtain
\begin{equation}\label{eq:5268}
\begin{split}
T_1^{**} \equiv & \: t_{\sigma^*_{M_{i+1}}} - t_{P_1}  =
\frac{\sqrt{2 \sigma_{M_{i}}}}{U-\Phi} \, \times \\
 & \sqrt{U(1-\beta_1) + \Phi(1-2\beta_1 +\beta_2)} \end{split}
\end{equation}
and
\begin{equation}\label{eq:5269}
T_2^{**} \equiv t_{\sigma^{**}_{M_{i+1}}} - t_{P_2} =
\frac{\sqrt{2 \sigma_{M_{i}}}}{U-\Phi} \sqrt{U(1-\beta_1) +
\Phi(1-\beta_2)},
\end{equation}
respectively. Also here one can show, cf. \eqref{eq:5268},
\eqref{eq:5269}, that
$$
T_1^{**} < T_2^{**}
$$
holds always due to $\beta_2 < \beta_1$.

For analyzing the total worst-case (i.e. largest) reaching time
between two consecutive extreme values $\sigma_{M_{i}}$ and
$\sigma_{M_{i+1}}$, one needs to distinguish between two pathes:
$P'' \rightarrowtail P_1 \rightarrowtail \sigma^*_{M_{i+1}}$ and
$P'' \rightarrowtail P_2 \rightarrowtail \sigma^{**}_{M_{i+1}}$,
cf. Figure \ref{fig:22}. In the first case, the overall maximal
reaching time is
\begin{equation} \label{eq:527}
\begin{split}
T_{1,\max,i} = & \:  T' + T^*_{1,\max} + T_1^{**} = \sqrt{2
\sigma_{M_{i}}} \, \Omega_1 =  \\
& \sqrt{2 \sigma_{M_{i}}} \, \Bigl( \bar{\Omega}_1 +
\tilde{\Omega}_1 \Bigr), \\
\hbox{with } \; \Omega_1 = &  \frac{1}{\Phi} \sqrt{(U+\Phi)(1-\beta_1)} \, + \\
& \frac{1}{U-\Phi} \sqrt{(U-\Phi)(1-\beta_1)} \, + \\
& \frac{1}{U-\Phi} \sqrt{U(1-\beta_1)+ \Phi(1-2\beta_1+\beta_2)} .
\end{split}
\end{equation}
In the second case, respectively, the overall maximal reaching
time is
\begin{equation} \label{eq:5275}
\begin{split}
T_{2,\max,i} = & \: T' + T^*_{2,\max} + T_2^{**} = \sqrt{2
\sigma_{M_{i}}} \, \Omega_2 = \\
& \sqrt{2 \sigma_{M_{i}}} \, \Bigl( \bar{\Omega}_2 +
\tilde{\Omega}_2 \Bigr), \\
\hbox{with } \; \Omega_2 = &  \frac{1}{\Phi} \sqrt{U(1-\beta_1) + \Phi(1-\beta_2)} \, + \\
& \frac{1}{U-\Phi} \sqrt{(U-\Phi)(1-\beta_1)} \, + \\
& \frac{1}{U-\Phi} \sqrt{U(1-\beta_1)+ \Phi(1-\beta_2)}.
\end{split}
\end{equation}
Note that in both above equations $\bar{\Omega}$ denote the terms
where the control is off, i.e. $u=0$, while $\tilde{\Omega}$
denote the terms where the control is on, i.e. $|u| = U$.

Once the convergence of the energy-saving sub-optimal SMC is
guaranteed, being analyzed in Section \ref{sec:3:sub:2}, the
overall accumulated worst-case reaching time results in
\begin{equation}\label{eq:59}
T_{\max,\sum} = \sqrt{2} \,  \max \{\Omega_1,  \Omega_2 \} \, \sum
\limits_{i=1}^{\infty} \sqrt{|\sigma_{M_{i}}|}.
\end{equation}
Because the convergence implies
\begin{equation}\label{eq:591}
|\sigma_{M_{i+1}}| = \eta |\sigma_{M_{i}}| \quad \hbox{with} \quad
0 < \eta < 1,
\end{equation}
the infinite sum in \eqref{eq:59} turns out into a geometric
series
\begin{equation}\label{eq:60}
\sum \limits_{i=1}^{\infty} \sqrt{|\sigma_{M_{i}}|} =
|\sigma_{M_{1}}|^{\frac{1}{2}} \sum \limits_{n=0}^{\infty}
\eta^{\frac{1}{2}n} = |\sigma_{M_{1}}|^{\frac{1}{2}} \frac{ 1
}{1-\eta^{\frac{1}{2}}},
\end{equation}
which is convergent. This allows concluding that the upper bound
of the finite time convergence is
\begin{equation}\label{eq:61}
T_c \leq t_{M_1} + \frac{\sqrt{2 } \, \max \{\Omega_1,  \Omega_2
\} }{1 - \sqrt{\eta}} \, \sqrt{|\sigma_{M_{1}}|}.
\end{equation}
Note that the initial reaching time of the first extreme value is
$t_{M_1} \leq |\dot{\sigma}(0)| / (U-\Phi)$, cf. \eqref{eq:9b}.
Furthermore, we note that depending on the branching of
trajectories, as shown above for the control-off phase, also a
maximum value $\eta = \max \{ \eta_1, \eta_2 \}$ must be taken
into account in \eqref{eq:61}.

For determining $\eta_1$ and $\eta_2$, we evaluate both subsequent
extreme values and obtain, based on \eqref{eq:10}, the equations
\begin{equation}\label{eq:63}
\begin{split}
\sigma^{*}_{M_{i+1}} = & \beta_2^{*} \sigma_{M_{i}} + 0.5 \frac{
-\dot{\sigma}_{P_1}^2}{U-\Phi} = \\
& \sigma_{M_{i}} \Biggl( \beta_2 - \frac{U(1-\beta_1) +
\Phi(1-2\beta_1+\beta_2)}{U-\Phi} \Biggr)
\end{split}
\end{equation}
and
\begin{equation}\label{eq:64}
\begin{split}
\sigma^{**}_{M_{i+1}} = & \beta_2^{**} \sigma_{M_{i}} + 0.5
\frac{- \dot{\sigma}_{P_2}^2}{U-\Phi} = \\
& \sigma_{M_{i}} \Biggl( \beta_2 - \frac{U(1-\beta_1) +
\Phi(1-\beta_2)}{U-\Phi} \Biggr),
\end{split}
\end{equation}
Comparing \eqref{eq:63} and \eqref{eq:64} with \eqref{eq:591},
results in
\begin{equation}\label{eq:65}
\eta_1 = \Biggl| \beta_2 - \frac{U(1-\beta_1) +
\Phi(1-2\beta_1+\beta_2)}{U-\Phi} \Biggr| < 1
\end{equation}
and
\begin{equation}\label{eq:66}
\eta_2 = \Biggl| \beta_2 - \frac{U(1-\beta_1) +
\Phi(1-\beta_2)}{U-\Phi} \Biggr| < 1.
\end{equation}

Following the same line of calculations as above, one can evaluate
the reaching and convergence time of the conventional sub-optimal
SMC, for which $\beta \equiv \beta_1$. Considering the $P'$ and
$P''$ peaking points results in
\begin{equation} \label{eq:528}
\begin{split}
\hat{T}_{\max,i} = & T' + \hat{T}^{**} = \sqrt{2 \sigma_{M_{i}}}
\, \hat{\Omega}, \\
\hbox{with } \; \hat{\Omega} = &  \frac{1}{U-\Phi}
\sqrt{(U+\Phi)(1-\beta_1)}\, + \\
& \frac{1}{U-\Phi} \sqrt{(U - \Phi)(1-\beta_1)}.
\end{split}
\end{equation}
Then, the upper bound of the finite time convergence of the
conventional sub-optimal SMC results in
\begin{equation}\label{eq:62}
\hat{T}_c < t_{M_1} + \frac{\sqrt{2 } \, \hat{\Omega}}{1 -
\sqrt{\hat{\eta}}} \, \sqrt{|\sigma_{M_{1}}|},
\end{equation}
where $0 < \hat{\eta} < 1$, cf. with \eqref{eq:591}. For
evaluating $\hat{\eta}$, consider the peaking point $P''$ and
obtain
\begin{equation}\label{eq:67}
\hat{\eta} = \Biggl| \beta_1 - \frac{(U+\Phi)(1-\beta_1)}{U-\Phi}
\Biggr| < 1.
\end{equation}
Note that both \eqref{eq:65} and \eqref{eq:66} reduce to
\eqref{eq:67} if assuming $\beta_2 = \beta_1$. Also worth noting
is that \eqref{eq:67} coincides with the convergence condition
\eqref{eq:18} if $\beta_2 = \beta_1$ and requires
\begin{equation}\label{eq:68}
\beta_1 > \frac{\Phi}{U}.
\end{equation}

We will denote $\Omega_1, \, \Omega_2, \, \hat{\Omega}$ as
\emph{amplification factor}s of the convergence time, and $\eta_1,
\, \eta_2, \, \hat{\eta}$ as \emph{contraction factor}s of the
convergence. Inspecting \eqref{eq:61} and \eqref{eq:62}, one can
recognize that the former are proportionally increasing the
convergence time while the latter can drastically slow down the
overall convergence once $\eta, \hat{\eta}$ become closer to one.

\subsection{Energy-saving control parametrization}
\label{sec:3:sub:4}

For the energy-saving sub-optimal SMC becomes effective, one needs
to ensure that its total energy consumption is less than that of
the conventional sub-optimal SMC, cf. condition (b) at the end of
Section \ref{sec:2}, while assuming the $U$, $\Phi$, and $\beta_1$
values are the same. In terms of the both control-on phases, i.e.
$u \neq 0$, the amplification factors to be taken into account are
\begin{equation}\label{eq:69}
\begin{split}
\tilde{\Omega}_1 = & \frac{1}{U-\Phi} \sqrt{(U-\Phi)(1-\beta_1)}\, + \\
& \frac{1}{U-\Phi} \sqrt{U(1-\beta_1)+ \Phi(1-2\beta_1+\beta_2)}
\end{split}
\end{equation}
and
\begin{equation}\label{eq:70}
\begin{split}
\tilde{\Omega}_2 = & \frac{1}{U-\Phi} \sqrt{(U-\Phi)(1-\beta_1)}\, + \\
& \frac{1}{U-\Phi} \sqrt{U(1-\beta_1)+ \Phi(1-\beta_2)},
\end{split}
\end{equation}
cf. \eqref{eq:527}, \eqref{eq:5275}. Taking out of consideration
the initial reaching phase in \eqref{eq:61} and \eqref{eq:62}, we
obtain two energy cost functions
\begin{equation}\label{eq:71}
J(\beta_1, \beta_2, U, \Phi) = \frac{ \max \{\tilde{\Omega}_1,
\tilde{\Omega}_2 \} }{1 - \max \{ \sqrt{\eta_1}, \sqrt{\eta_2} \}
}
\end{equation}
and
\begin{equation}\label{eq:72}
\hat{J}(\beta_1, U, \Phi) = \frac{ \hat{\Omega} }{1 -
\sqrt{\hat{\eta}} },
\end{equation}
for the energy-saving and conventional sub-optimal SMC,
respectively. In order to find an energy-saving pair of the
$(\beta_1,\beta_2)$ threshold parameters for with some fixed
$U,\Phi$ process values satisfying $U>\Phi$, we formulate the
\begin{equation}\label{eq:72a}
\underset{\beta_1,\beta_2} {\min} \, \Bigl[  J(\beta_1, \beta_2,
U, \Phi) - \hat{J}(\beta_1, U, \Phi) \Bigr]
\end{equation}
minimization problem under the hard constraint
\begin{equation}\label{eq:72b}
J(\beta_1, \beta_2, U, \Phi) - \hat{J}(\beta_1, U, \Phi) < 0.
\end{equation}
In addition, in order to avoid $\beta_1 \rightarrow 1$, which can
lead to $\hat{T}_c \rightarrow \infty$, see \cite{athans1964} for
details, we specify the upper bound
\begin{equation}\label{eq:72c}
\hat{J}(\beta_1, U, \Phi) < \hat{J}_{\max}.
\end{equation}
This implies restriction on the convergence time $\hat{T}_c$ and,
this way, forces a reasonable $\beta_1 < 1$ in course of solving
the minimization problem \eqref{eq:72a}.
\begin{figure}[!h]
\centering
\includegraphics[width=0.98\columnwidth]{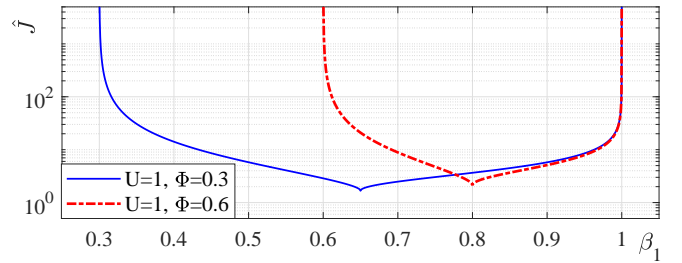}
\caption{Convergence time cost function $\hat{J}$ in dependency of
$\beta_1$.} \label{fig:31}
\end{figure}
The convergence time cost function $\hat{J}$ is exemplary
visualized in Figure \ref{fig:31} in dependency of $\beta_1$, for
two cases assuming $U=1$ and $\Phi=\{0.3, \, 0.6 \}$. Note that
$\hat{J} \rightarrow \infty$ for both boundary values $\beta_1
\rightarrow \Phi/U$, cf. \eqref{eq:68}, and $\beta_1 \rightarrow
1$.

\begin{figure}[!h]
\centering
\includegraphics[width=0.98\columnwidth]{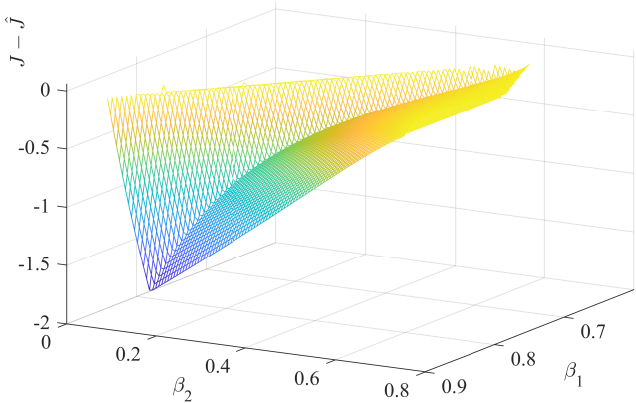}
\caption{Constrained objective function $(J-\hat{J})$ of the
$\beta_1$, $\beta_2$ parameters, for exemplary perturbation to
control ratio $\Phi/U=0.3$.} \label{fig:32}
\end{figure}

Note that the hard constraint \eqref{eq:72b} is required to ensure
that the total control-on time of the energy-saving sub-optimal
SMC is shorter than that of the conventional one. Moreover, one
can recognize that \eqref{eq:72a} ensures the maximal possible
energy-saving, due to the control-off phases throughout the
convergence, comparing to the conventional sub-optimal SMC, both
having the same $\beta_1$ value. That means for any fixed
$\beta_1$ complying with \eqref{eq:72c}, there is an optimal
$\beta_2$ counterpart. The latter represents the solution of the
constrained minimization problem \eqref{eq:72a}. The results of
the constrained optimization are exemplary visualized in Figure
\ref{fig:32}, for the perturbation to control ratio equal to
$\Phi/U=0.3$.
It should be noted that the formulated optimization problem does
not take explicitly into account those transient phases where the
control is off, cf. \eqref{eq:527}, \eqref{eq:5275} and
\eqref{eq:69}, \eqref{eq:70}. However, this does not imply any
issues since the finite time convergence of the energy-saving
sub-optimal SMC is upper bounded by \eqref{eq:61}.

\section{Chattering due to actuator dynamics}
\label{sec:4}

If the system under control is subject to an additional
(parasitic) actuator dynamics, the second-order sliding mode
experiences the residual steady-state oscillations, also known as
\emph{chattering}, see e.g. \cite{bartolini1998,boiko2007}. For
the first-order actuator dynamics with a time constant $\mu > 0$,
the control variable $u(t)$ in \eqref{eq:2} must be substituted by
a new control variable $v(t)$, which is the solution of
$$
\mu \dot{v}(t) + v(t) = u(t).
$$

\begin{figure}[!h]
\centering
\includegraphics[width=0.65\columnwidth]{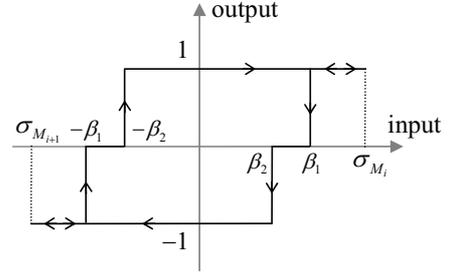}
\caption{Hysteresis relay representation of the control
\eqref{eq:9a} at residual steady-state oscillations due to the
parasitic actuator dynamics.} \label{fig:41}
\end{figure}

Assume that the stable residual steady-state oscillations of the
system \eqref{eq:2} with the new input channel $v$ and control
\eqref{eq:9a} are established. Then, the input-output map of the
energy-saving sub-optimal SMC \eqref{eq:9a} takes the form of a
three-state hysteresis relay as shown in Figure \ref{fig:41}. For
steady-state oscillations, the switching thresholds keep the
constant values and are symmetric with respect to zero. Note that
the amplitude of oscillations is determined by the cyclic extreme
value $\sigma_{M_i} \equiv \sigma_A$. Furthermore, we note that if
$\beta_1=\beta_2$, the input-output map reduces to a standard
two-state hysteresis relay, as it was used for harmonic balance
analysis of the original sub-optimal control, cf.
\cite[chapter~5.2]{boiko2009}.

\begin{figure}[!h]
\centering
\includegraphics[width=0.7\columnwidth]{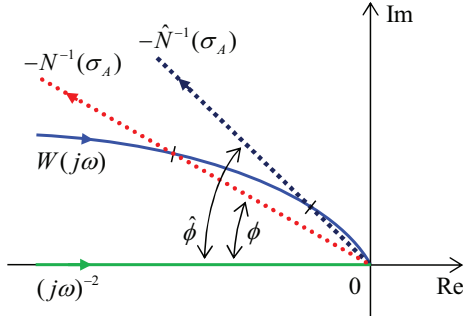}
\caption{Negative reciprocal of DF of the conventional and
energy-saving sub-optimal SMC crossing the Nyquist plot of
$W(j\omega)$.} \label{fig:42}
\end{figure}

Before examining the residual oscillating behavior of the
energy-saving sub-optimal SMC, recall the principle procedure of
the describing function (DF) analysis, see e.g. \cite{Atherton75},
here in context of the second-order SMC and, in particular, when
using the sub-optimal SMC. For the given DF of the SMC in
feedback, denoted by $N(\sigma_A)$, the harmonic balance equation
\begin{equation}\label{eq:73}
N(\sigma_A) \, W(j\omega) + 1 = 0
\end{equation}
must have a real solution, in terms of the amplitude $\sigma_A$
and angular frequency $\omega$, for the residual steady-state
oscillations to exist. Considering the double-integrator plant
augmented by the first-order actuator dynamics, we write
\begin{equation}\label{eq:74}
W(j\omega) = \frac{1}{s^2(\mu s + 1)} \Biggr|_{s=j\omega} =
\frac{-1}{\omega^2 + j \cdot \omega^3 \mu}.
\end{equation}
For the conventional sub-optimal SMC, parameterized by $U$ and
$\beta_1$, the DF is given by, cf. \cite{boiko2009},
\begin{equation}\label{eq:75}
\hat{N}(\sigma_A) = \frac{4U}{\pi \sigma_A}
\Bigl(\sqrt{1-\beta_1^2} + j \beta_1 \Bigr).
\end{equation}
Note that also for the generalized sub-optimal SMC, the DF
analysis was provided in \cite{boiko2006}, cf. with \eqref{eq:75}.
The negative reciprocal of \eqref{eq:75} yields
\begin{equation}\label{eq:75b}
-\hat{N}(\sigma_A)^{-1} = \frac{\pi \sigma_A}{4U}
\Bigl(-\sqrt{1-\beta_1^2} + j \beta_1 \Bigr).
\end{equation}
One can show that the locus of \eqref{eq:75b}, that is its
graphical interpretation in the harmonic balance \eqref{eq:73}, is
a straight line starting in the origin and forming a clockwise
angle $\hat{\phi}$, see Figure \ref{fig:42}. The latter can be
directly calculated as
\begin{equation}\label{eq:75c}
\hat{\phi} = \arctan \Biggl(-\frac{\beta_1}{\sqrt{1-\beta_1^2}}
\Biggr).
\end{equation}
One can recognize that larger $\beta_1$-parameter values lead to a
larger $|\hat{\phi}|$ and, therefore, shirt the locus intersection
towards higher frequencies and lower amplitudes.

Since the energy-saving sub-optimal SMC constitutes a linear
combination of two negative hysteresis relays, cf. \eqref{eq:9a},
its DF can be represented as a summation of two DFs of the type
\eqref{eq:75}, this way resulting in
\begin{equation}\label{eq:76}
N(\sigma_A) = \frac{2U}{\pi \sigma_A} \Bigl(\sqrt{1-\beta_1^2} +
\sqrt{1-\beta_2^2} + j (\beta_1 + \beta_2) \Bigr).
\end{equation}
After evaluating the negative reciprocal of \eqref{eq:76}, similar
to \eqref{eq:75b}, we obtain
\begin{equation}\label{eq:77}
\phi = \arctan
\Biggl(-\frac{\beta_1+\beta_2}{\sqrt{1-\beta_1^2}+\sqrt{1-\beta_2^2}}
\Biggr).
\end{equation}
In Figure \ref{fig:42}, the negative reciprocal of both DFs are
depicted, once of the conventional sub-optimal SMC with DF
\eqref{eq:75} and once of the energy-saving sub-optimal SMC with
DF \eqref{eq:76}, together with the Nyquist plots of $W(j\omega)$
and double integrator $(j\omega)^{-2}$. The energy-saving
sub-optimal SMC has always a lower angle $|\phi|$ for the same
fixed $\beta_1$ value, while for $(\beta_1-\beta_2) \rightarrow 0$
the angle $\phi \rightarrow \hat{\phi}$. One can also recognize
that unlike the double integrator, the Nyquist plot of
$W(j\omega)$ has always an intersection point with the negative
reciprocals of DFs. This implies an existence of $(\sigma_A,
\omega)$ solutions of the harmonic balance equation \eqref{eq:73}.
The corresponding steady-state oscillations have a higher
amplitude and lower frequency for the energy-saving sub-optimal
SMC compared to those of the conventional sub-optimal SMC.

Equating the phase angle of \eqref{eq:74} and \eqref{eq:77}
$$
\frac{\mathrm{Im} \{ W \}}{\mathrm{Re} \{ W \}} =
-\frac{\beta_1+\beta_2}{\sqrt{1-\beta_1^2}+\sqrt{1-\beta_2^2}}
$$
and solving it with respect to $\omega$ results in
\begin{equation}\label{eq:78}
\omega_c = \mu^{-1} \frac{\beta_1+\beta_2}{
\sqrt{1-\beta_1^2}+\sqrt{1-\beta_2^2}}
\end{equation}
of the chattering. The corresponding amplitude of the steady-state
oscillations is then obtained as
\begin{equation}\label{eq:79}
\sigma_A = \frac{\sqrt{\mu^2 \omega_c^2 + 1}}{\omega_c^2 \bigl(
\mu^2 \omega_c^2 + 1 \bigr) }.
\end{equation}

\section{Numerical examples}
\label{sec:5}

This section demonstrates the numerical examples which
characterize the main properties of the energy-saving sub-optimal
SMC formulated and analyzed above. Implemented is a
double-integrator system \eqref{eq:2} with $g=1$ and the
energy-saving sub-optimal SMC \eqref{eq:9a}. For simulating the
conventional sub-optimal SMC, the threshold value is set
$\beta_2=\beta_1$, cf. section \ref{sec:2}. All simulation results
are obtained by using the simple first-order Euler-type numerical
solver with the fixed-step size of 0.001 sec.

\subsection{Convergence behavior}
\label{sec:5:sub:1}

\begin{figure}[!h]
\centering
\includegraphics[width=0.98\columnwidth]{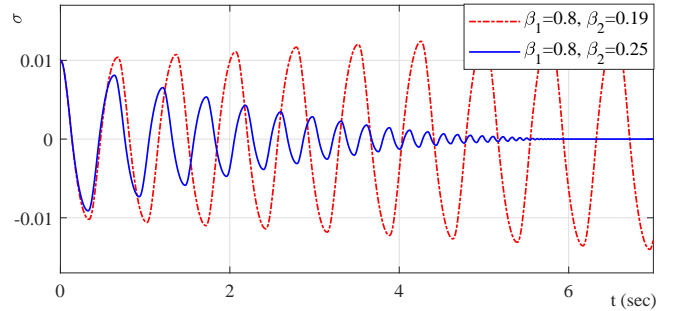}
\caption{Exemplary response of the controlled sliding variable for
$\beta_1, \beta_2$ once fulfilling (in blue) and once violating
(in red) the parametric convergence condition \eqref{eq:18};
$U=1$, $\Phi=0.5$, $f = 0.5 \, \mathrm{sign}(\dot{\sigma})$.}
\label{fig:51}
\end{figure}

The convergence behavior of the energy-saving sub-optimal SMC is
exemplary shown for examining the parametric condition
\eqref{eq:18}. The constraining parameters are chosen so that
$\Phi/U=0.5$ with $U=1$ and $\mu=0$, while for a worst-case
scenario of the bounded perturbation $f = 0.5 \,
\mathrm{sign}(\dot{\sigma})$ is assumed. Note that the assumed $f$
constitutes an always 'co-acting' perturbation signal of the
maximal amplitude $\Phi$, a worst-case scenario that challenges
the contraction during each reaching cycle and, at large, the
total convergence to $\sigma=\dot{\sigma}=0$ equilibrium. Figure
\ref{fig:51} discloses the converging behavior for $\beta_1 =
0.8$, $\beta_2 = 0.25$ versus the diverging behavior for $\beta_1
= 0.8$, $\beta_2 = 0.19$. Note that here, the threshold values
distance $\beta_1-\beta_2 > 0.6$ is required for \eqref{eq:18} is
fulfilled.

\subsection{Energy-saving versus conventional sub-optimal SMC}
\label{sec:5:sub:2}

The energy-saving sub-optimal control behavior is analyzed in
comparison to the conventional sub-optimal one, with regard to the
recorded convergence time $T_c$ and overall energy consumption
$$
\int_{0}^{T_c} |u(t)| dt.
$$
In order to determine the numerical convergence instant, the $T_c$
value is obtained once the Euclidian norm of the state vector
$\bigl(\sigma(t), \dot{\sigma}(t)\bigr)$ decreases below some
assigned low residual constant, i.e. $\bigl\| (\sigma,
\dot{\sigma}) \bigr\|(t) < 0.004$. The residual constant is
determined out from the numerical simulations, with respect to the
solver type and step-size and the resulted numerical pattern
around $\sigma=0$ after the principal dynamics of $\sigma(t)$ has
converged. For the assigned $\Phi/U=0.3$ ratio with $U=1$ and
$\mu=0$, three configurations of the perturbing signal are
considered: $f=0$, $f=0.3 \, \mathrm{sign}(\dot{\sigma})$, $f=-0.3
\, \mathrm{sign}(\dot{\sigma})$. Note that the second one is
'co-acting' to the control signal. This has, however, an adverse
impact on contraction of the state trajectories, cf. Figure
\ref{fig:22}. The third one acts, on the contrary, as a maximal
possible system damping and can negatively affect the control-off
phases in terms of the reaching time, cf. section
\ref{sec:3:sub:3}. Two parameter sets, obtained from the
constrained optimization, $[\beta_1=0.7, \beta_2=0.55]$ and
$[\beta_1=0.83, \beta_2=0.32]$ cf. with Figure \ref{fig:32}, are
used. All simulation results are summarized in Table \ref{tab:0}.
Recall that for the conventional sub-optimal SMC the $\beta_2 =
\beta_1$ is always assigned.
\begin{table}[!h]
  \renewcommand{\arraystretch}{1.4}
  \caption{Convergence time and energy consumption}
  \label{tab:0}
  \scriptsize
  \begin{center}
  \begin{tabular} {|p{2.1cm}||p{1.1cm}|p{1.1cm}||p{1.1cm}|p{1.1cm}|}
  \hline
  \multicolumn{1}{| c ||}{} & \multicolumn{2}{c||} {$\beta_1=0.7$, $\beta_2=0.55$} & \multicolumn{2}{c|} {$\beta_1=0.83$, $\beta_2=0.32$} \\
  \cline{2-5}
  Configuration          & $T_c$   &  $\int |u|dt$   & $T_c$  & $\int |u|dt$ \\
  \hline \hline
  conventional \newline $f=0$                                    & 0.420      &  0.391         & 0.643                   & 0.580       \\
  \cline{1-5}
  energy-saving \newline $f=0$                                   & 0.345      &  0.292         & 0.332                   & 0.176       \\
  \hline \hline
  conventional \newline  $f=0.3 \, \mathrm{sign}(\dot{\sigma})$  & 0.342      &  0.285         & 0.552                   & 0.448       \\
  \cline{1-5}
  energy-saving \newline $f=0.3 \, \mathrm{sign}(\dot{\sigma})$  & 0.288      &  0.247         & 0.435                   & 0.278       \\
  \hline \hline
  conventional \newline $f=-0.3 \, \mathrm{sign}(\dot{\sigma})$  & 0.517      &  0.441         & 0.645                   & 0.548       \\
  \cline{1-5}
  energy-saving \newline $f=-0.3 \, \mathrm{sign}(\dot{\sigma})$ & 0.442      &  0.325         & \cellcolor{red!25}3.718 & 0.155       \\
  \hline \hline
  \end{tabular}
  \end{center}
\end{table}
One can recognize a superior performance of the energy-saving
sub-optimal SMC in all above given cases, except the last row in
Table \ref{tab:0} (see light red shadowed). Here the energy saving
is largely achieved, but the convergence time is significantly
increased comparing to the conventional sub-optimal SMC. This is
related to a relatively slowly 'drifting' trajectory between the
lying apart $\beta_1=0.83$ and $\beta_2=0.32$, being driven by the
$f=-0.3 \, \mathrm{sign}(\dot{\sigma})$ perturbation. The
situation with $T_c$ improves, still guaranteing for energy
saving, when both threshold parameters are lying closer to each
other like in the shown case $\beta_1=0.7$ and $\beta_2=0.55$, see
Table \ref{tab:0}.

\begin{figure}[!h]
\centering
\includegraphics[width=0.49\columnwidth]{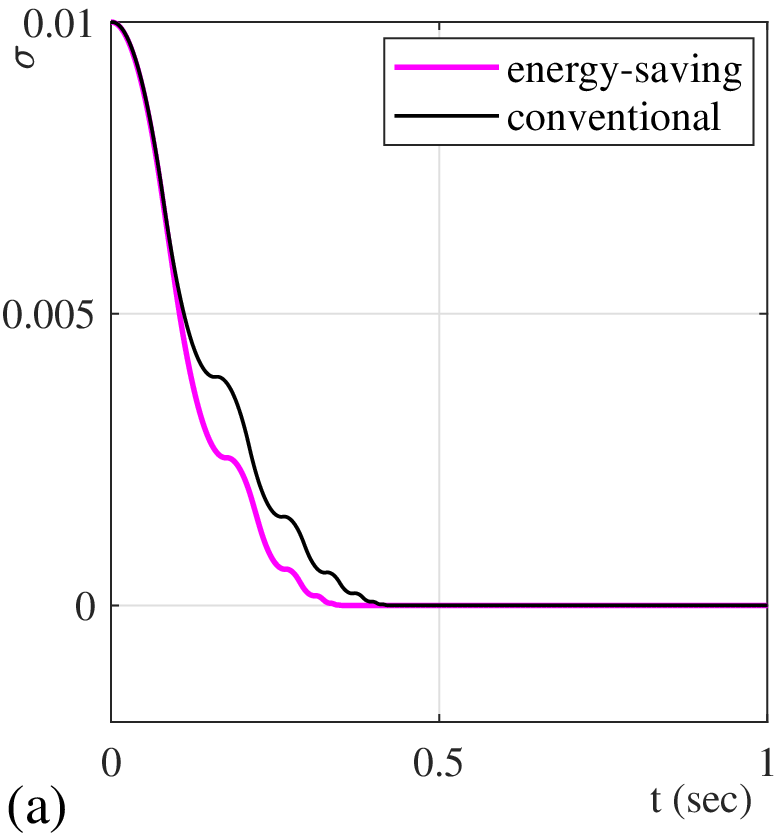}
\includegraphics[width=0.49\columnwidth]{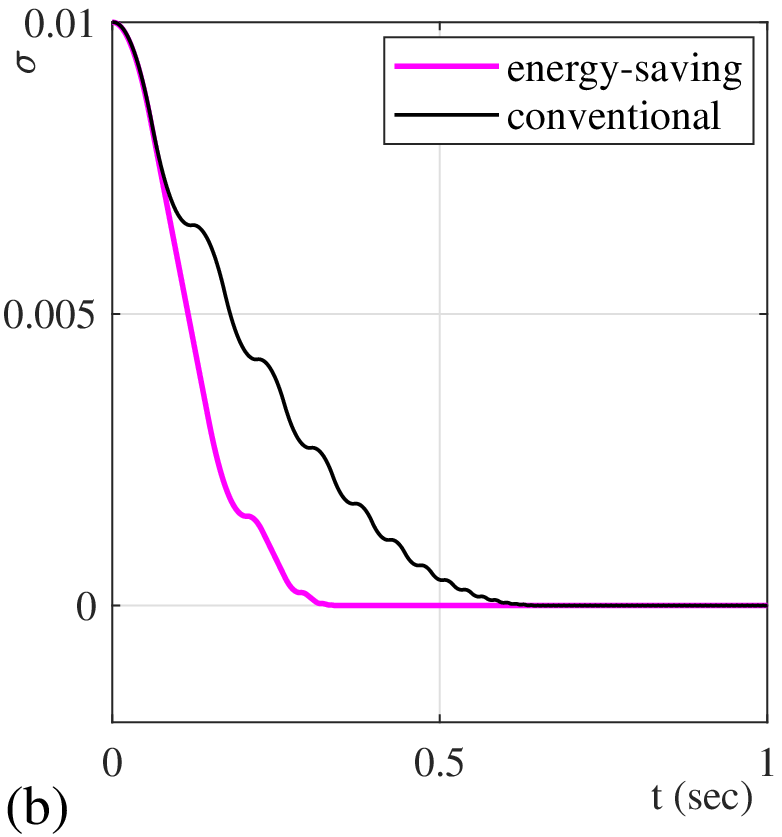}
\caption{Convergence of $\sigma$ with $f=0$ for the parameters
$[\beta_1=0.7, \beta_2=0.55]$ in (a) and $[\beta_1=0.83,
\beta_2=0.32]$ in (b).} \label{fig:521}
\end{figure}

\begin{figure}[!h]
\centering
\includegraphics[width=0.49\columnwidth]{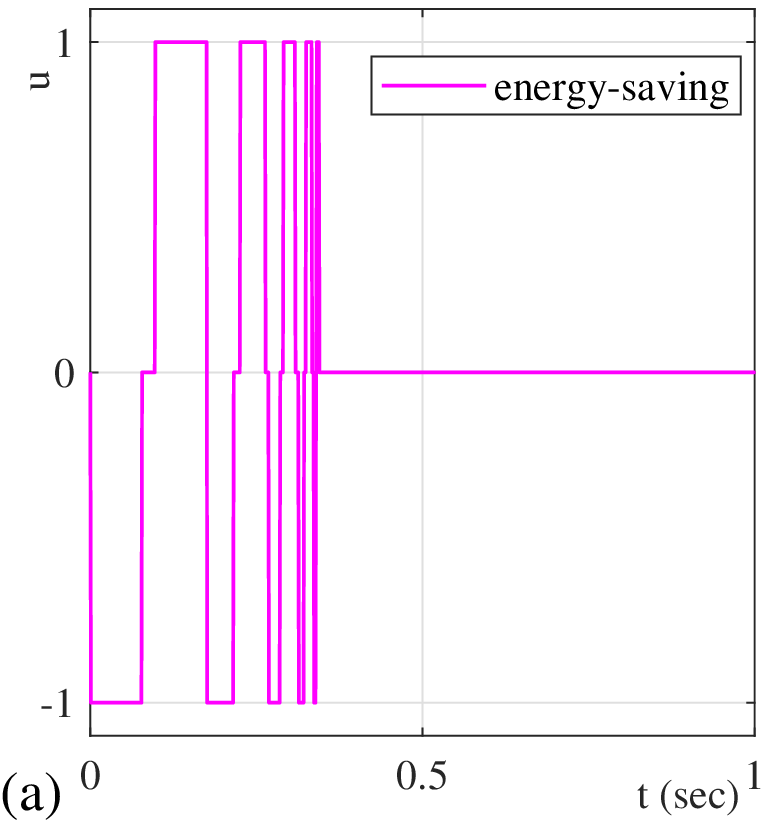}
\includegraphics[width=0.49\columnwidth]{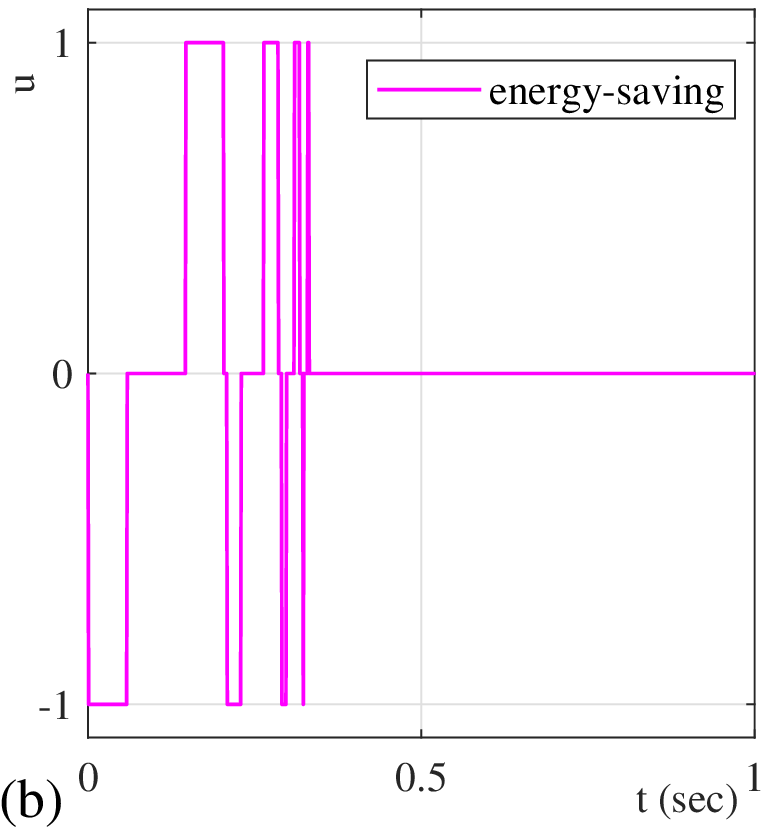}
\includegraphics[width=0.49\columnwidth]{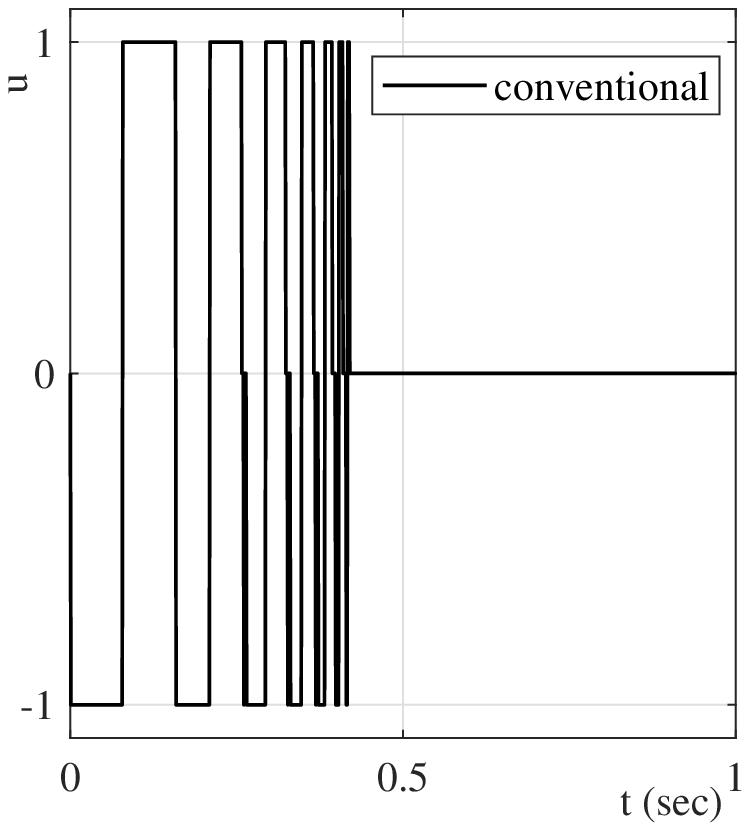}
\includegraphics[width=0.49\columnwidth]{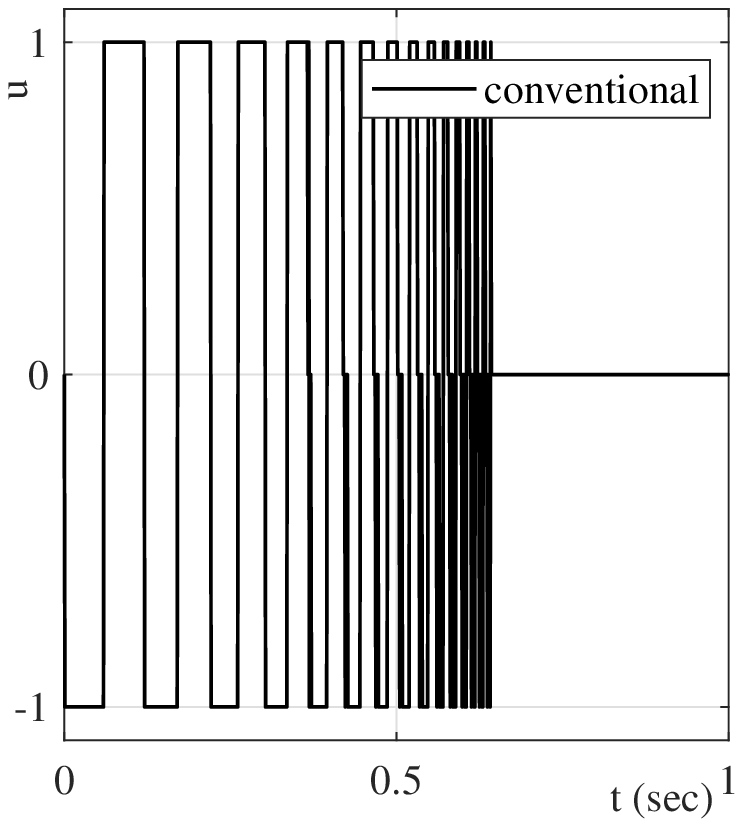}
\caption{Control $u$ with $f=0$ for the parameters $[\beta_1=0.7,
\beta_2=0.55]$ in (a) on the left, and $[\beta_1=0.83,
\beta_2=0.32]$ in (b) on the right.} \label{fig:522}
\end{figure}

\begin{figure}[!h]
\centering
\includegraphics[width=0.49\columnwidth]{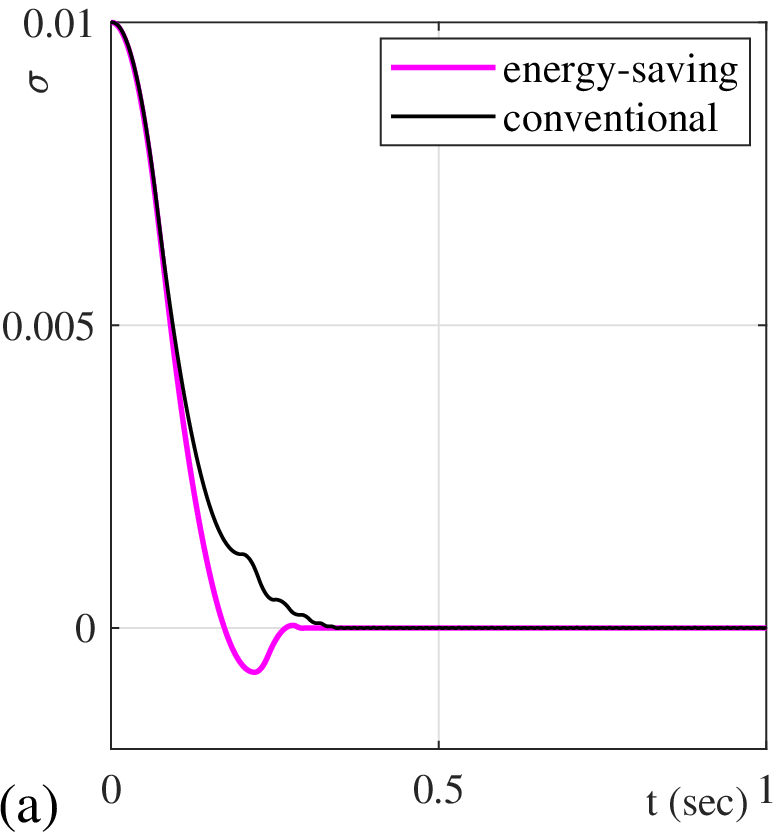}
\includegraphics[width=0.49\columnwidth]{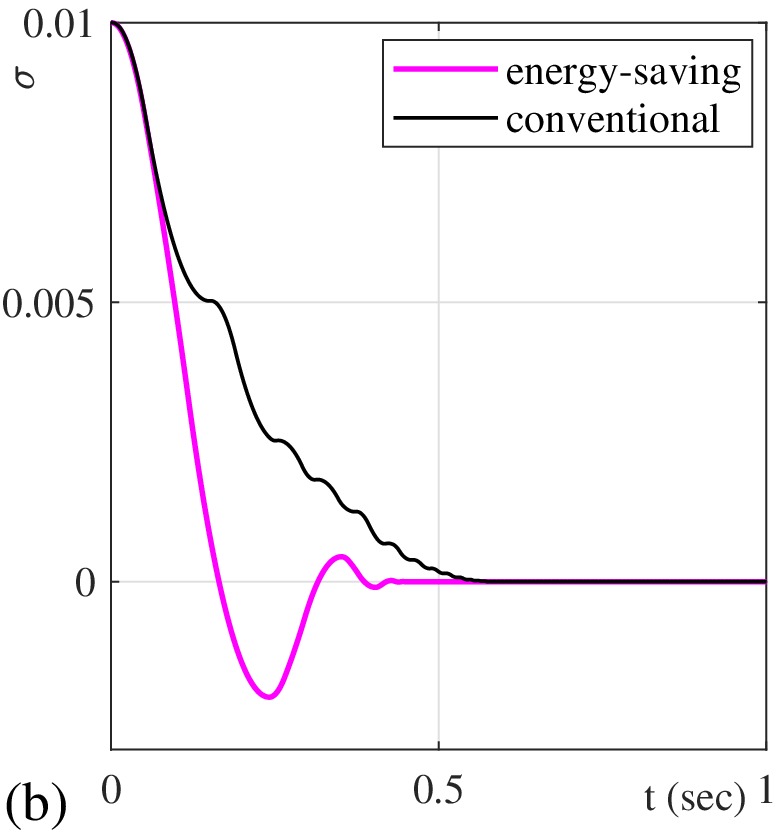}
\caption{Convergence of $\sigma$ with $f=0.3 \,
\mathrm{sign}(\dot{\sigma})$ for the parameters $[\beta_1=0.7,
\beta_2=0.55]$ in (a) and $[\beta_1=0.83, \beta_2=0.32]$ in (b).}
\label{fig:523}
\end{figure}

For better interpreting the listed results, some of the above
configurations are shown below in the plots of $\sigma(t)$ and
$u(t)$. All configurations with $f=0$ are depicted in Figures
\ref{fig:521} and \ref{fig:522}.
One can see that the energy-saving sub-optimal SMC converges
faster, correspondingly stops to switch earlier than its
conventional counterpart. One can also notice that both
controllers here disclose a monotonic convergence without
twisting, cf. section \ref{sec:2}.
The convergence in case of $f=0.3 \, \mathrm{sign}(\dot{\sigma})$
is also exemplary shown in Figure \ref{fig:523}. Here, the
conventional sub-optimal SMC has a monotonic convergence, in
accord with \eqref{eq:9}. On the contrary, the energy-saving
sub-optimal SMC changes into twisting mode here, because it allows
for control phases driven only by the perturbation quantity.

\subsection{Residual steady-state oscillations}
\label{sec:5:sub:3}

\begin{figure}[!h]
\centering
\includegraphics[width=0.95\columnwidth]{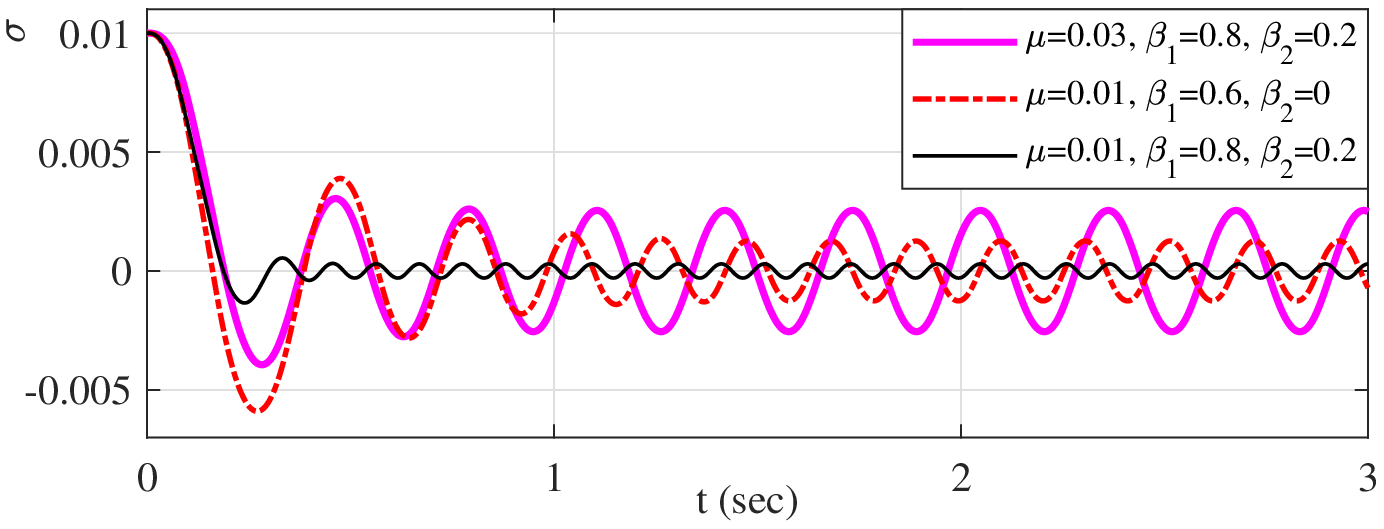}
\includegraphics[width=0.95\columnwidth]{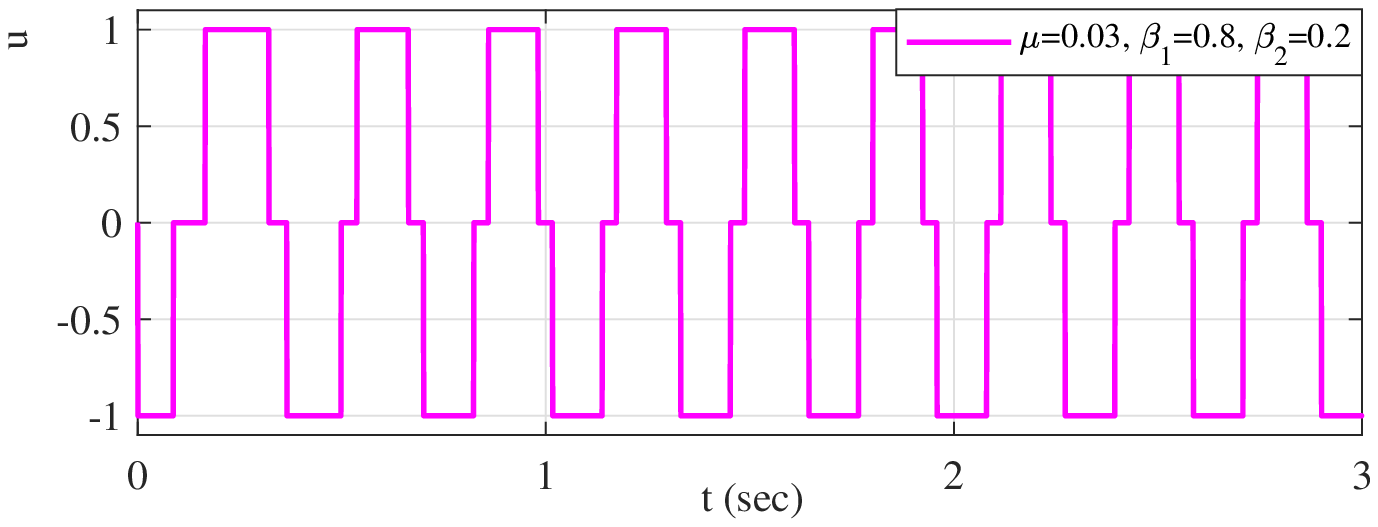}
\caption{Exemplary response of the residual steady-state
oscillations, sliding variable (above) and control value (below).}
\label{fig:53}
\end{figure}

For evaluating the residual oscillations (i.e. chattering) owing
to a parasitic actuator dynamics, cf. section \ref{sec:4},
different parameter sets with $\mu \neq 0$ are compared. Three
parameter configurations used are shown in Table \ref{tab:1}. The
amplitude and frequency of the residual steady-state oscillations
are computed based on the harmonic balance analysis provided in
section \ref{sec:4}, i.e. using \eqref{eq:78}, \eqref{eq:79}, and
compared with those out from the numerical simulations. Certain,
yet acceptable, error of the harmonic balance analysis is not
surprising, since DF solely approximates the first main harmonic,
cf. with results in \cite{boiko2006}.
\begin{table}[!h]
  \renewcommand{\arraystretch}{1.4}
  \caption{Harmonic balance analysis versus simulation results}
  \label{tab:1}
  \scriptsize
  \begin{center}
  \begin{tabular} {|p{2cm}||p{1.1cm}|p{1.1cm}||p{1.1cm}|p{1.1cm}|}
  \hline
  \multicolumn{1}{| c ||}{}                  & \multicolumn{2}{c||} {Harmonic balance} & \multicolumn{2}{c|} {Numerical simulation} \\
  \cline{2-5}
  Parameters                                 & $\sigma_A$   &  $\omega_c$   & $\sigma_A$  & $\omega_c$ \\
  \hline \hline
  $\mu=0.03$, $\beta_1=0.8$, $\beta_2=0.2$   & 0.0019       &  21.1         & 0.0025      & 20.0       \\
  \cline{1-5}
  $\mu=0.01$, $\beta_1=0.6$, $\beta_2=0.0$   & 0.00085      &  33.3         & 0.0012      & 30.2       \\
  \cline{1-5}
  $\mu=0.01$, $\beta_1=0.8$, $\beta_2=0.2$   & 0.00021      &  63.3         & 0.00029     & 59.3       \\
  \hline
  \end{tabular}
  \end{center}
\end{table}
Figure \ref{fig:53} exemplifies the appearance of residual
steady-state oscillations when the actuator dynamics is in place.
The time response of the sliding variable $\sigma(t)$ is shown
above for all three parameter sets, cf. Table \ref{tab:1}. The
control value $u(t)$, i.e. of the control \eqref{eq:9a}, is
exemplary shown below for the parameter set $\mu=0.03$,
$\beta_1=0.8$, $\beta_2=0.2$.
Table \ref{tab:1} and Figure \ref{fig:53} confirm the chattering
analysis of the energy-saving sub-optimal SMC, developed in
section \ref{sec:4}, and disclose how the amplitude and frequency
change depending on the control and actuator parameters $\beta_1$,
$\beta_2$, and $\mu$.

\section{Conclusions}
\label{sec:6}

The problem of saving energy by the control-off phases in the
second-order sliding modes was addressed. More specifically, the
conventional sub-optimal SMC, see \cite{bartolini2003} for survey,
was extended for a class of the systems \eqref{eq:2} with $g=1$
and $u \in \{-U,\, 0,\, U \}$. This way, a well parameterizable
$u=0$ phase appears within each reaching cycle, during the entire
finite-time convergence. The latter was estimated for the
worst-case scenario of an upper-bounded perturbation occurring in
each control phase. The analyzed contraction and maximal reaching-
as well as total convergence-time served to formulate a
constrained minimization problem for an energy cost function. Two
switching threshold parameters represent the solution of a complex
nonlinear optimization, while the perturbation-boundary to
control-authority ratio, i.e. $\Phi/U$, is crucial for convergence
and energy-saving conditions. Also the impact of an additional
parasitic actuator dynamics, typical for real sliding modes in
terms of chattering, was analyzed based on the harmonic balance.
The overall analysis performed and the dedicated numerical
evaluations argue in favor of using the proposed energy-saving
sub-optimal SMC for $\Phi/U < 0.35$. In applications, a resettable
control-off mode can also be added for $\bigl\| (\sigma,
\dot{\sigma}) \bigr\|(t) < \mathrm{const}$, this way ensuring
energy saving after convergence and a reactivation in case of
perturbations in the equilibrium.

\bibliographystyle{elsarticle-num}
\bibliography{references}

\end{document}